\documentclass{amsart}

\usepackage{amssymb}

\usepackage{url}
\newcommand{\R}{\mathbb{R}}

\newcommand{\C}{\mathbb{C}}
\newcommand{\Z}{\mathbb{Z}}

\newcommand{\N}{\Z_{>0}}
\newcommand{\bigoh}{{\mathcal{O}}}
\newcommand{\e}{{\rm e}}
\newcommand{\atan}{{\rm atan}}
\newcommand{\D}{{\rm d}}
\newcommand{\sinc}{{\rm sinc}}
\newcommand{\erf}{{\rm erf}}

\newtheorem{theorem}{Theorem}[section]
\newtheorem{lemma}[theorem]{Lemma}

\theoremstyle{definition}
\newtheorem{definition}[theorem]{Definition}

\theoremstyle{remark}

\numberwithin{equation}{section}

\begin{document}

\title{Numerical Computations Concerning the GRH}


%    author one information
\author{David J. Platt}
\address{Heilbronn Institute for Mathematical Research, University of Bristol, University Walk, Bristol, BS8 1TW, United Kingdom}
%\curraddr{}
\email{dave.platt@bris.ac.uk}
\thanks{These algorithms were developed during the author's Doctoral research under the patient supervision of Dr. Andrew Booker. The author in also indebted to the technical staff and management of ACRC at the University of Bristol, MesoPSL (affiliated with the Observatoire de Paris and Paris Sciences et Lettres), Direction des Syst\`{e}mes d'Information at Universit\'e Paris VI/VII (Pierre et Marie Curie), France Grilles (French National Grid
Infrastructure, DIRAC instance), PSMN at Universit\'e de Lyon 1 and PlaFRIM at Universit\'e de Bordeaux 1 for their invaluable support.}

%    author two information
%\author{}
%\address{}
%\curraddr{}
%\email{}
%\thanks{}

%    \subjclass is required.

\date{}

%\dedicatory{}

%    Abstract is required.
\begin{abstract}
We describe two new algorithms for the efficient and rigorous computation of
Dirichlet L-functions and their use to verify the Generalised
Riemann Hypothesis for all such L-functions associated with primitive characters of
modulus $q\leq 400,000$. For even $q$, we check to height
\begin{equation*}
t_0=\textrm{max}\left(\frac{10^8}{q},\frac{7.5\cdot 10^7}{q}+200\right)
\end{equation*}
 and for odd $q$ to height
\begin{equation*}
t_0=\textrm{max}\left(\frac{10^8}{q},\frac{3.75\cdot 10^7}{q}+200\right).
\end{equation*}
\end{abstract}

\maketitle

\section{Introduction}

For a given modulus $q\in\N$ we define the Dirichlet characters $\chi:\Z\to\C$ axiomatically
as follows:-
\begin{itemize}
\item $\chi(n)=0$ iff $(n,q)\neq 1$,
\item $\chi(mn)=\chi(m)\chi(n)$ and
\item $\chi(n+q)=\chi(n)$.
\end{itemize}
There are $\varphi(q)$ distinct characters of modulus $q$, where $\varphi$ is
Euler's totient function. The character $\chi(n)=1$ for all $n$ co-prime to $q$
is known as the principal character. A character $\chi$ of modulus $q$ is
primitive if and only if
for all $d$ dividing $q$ with $0<d<q$ there exists an integer $a\equiv 1\mod
d$ with $(a,q)=1$ and $\chi(a)\neq 1$. \cite{Apostol1976}. Finally, we define the parity of a character by
\begin{equation*}
a_\chi:=\frac{1-\chi(-1)}{2}.
\end{equation*}

The Dirichlet L-function of modulus $q$ associated with a character $\chi$ is
defined for  $\Re s>1$ by
\begin{equation*}
L_\chi(s)=\sum\limits_{n=1}^\infty \chi(n)n^{-s}
\end{equation*}
and with analytic continuation to $\C$ excepting (in the case of principal
characters) a simple pole at $s=1$.

Given $\epsilon_\chi$ such that $|\epsilon_\chi|=1$, we form the completed L-function via
\begin{equation*}
\Lambda_\chi(t):=\epsilon_\chi\left(\frac{q}{\pi}\right)^\frac{it}{2}\Gamma\left(\frac{\frac{1}{2}+a_\chi+it}{2}\right)\exp\left(\frac{\pi t}{4}\right)L_\chi\left(\frac{1}{2}+it\right).
\end{equation*}
For suitably chosen $\epsilon_\chi$, $\Lambda_\chi$ is real valued and has the same zeros as $L_\chi\left(\frac{1}{2}+it\right)$. The exponential factor is introduced (for computational expedience) to counteract the decay of the gamma function as $t$ increases.

The case $q=1$ we have only the principal character leading to a single L-function, namely Riemann's zeta function. Riemann's
guess that all zeros of this function with real part in $[0,1]$ lie on the
$1/2$ line is the Riemann Hypothesis (RH). Extensive calculations have been
undertaken to test RH to ever increasing heights, with Gourdon having checked
$10^{13}$ zeros \cite{Gourdon2010} using an algorithm first described by
Odlyzko and Sch{\"o}nhage \cite{Odlyzko1988}.

In contrast, the equivalent hypothesis for Dirichlet L-functions of primitive
character, which we will refer to as the
Generalised Riemann Hypothesis (GRH), has received less
attention. The last significant rigorous computation was that by Rumely \cite{Rumely1993}
who confirmed that the GRH holds for primitive L-functions modulus $q\leq 13$ to
height $10,000$ and various other moduli to height $2,500$.\footnote{Rumley
  refers the Extended Riemann Hypothesis, rather than the GRH. This former term is
  now more often used to describe the hypothesis as related to zeros of the
  Dedekind Zeta functions.}
The largest modulus tested was $q=432$ and in total about $10^7$ zeros were
examined. We note that Rumely went on to isolate these zeros with some
precision and to generate statistics on their locations, but in terms of
simply the number of zeros confirmed to lie on the $1/2$ line, there remained a factor of $10^6$ in
favour of zeta. If this weren't motivation enough, recent advances in the application of the Circle Method held out the tantalising prospect that ternary Goldbach might succumb to a combined numerical and analytic assault.

We will describe a computation using two new algorithms and exploiting improvements in hardware in the $20$ years since Rumely's paper that extend his result by about $6$ orders of magnitude in terms of the number of zeros checked. Furthermore, the combination of moduli and heights checked is more than sufficient to support Helfgott's proof of ternary Goldbach \cite{Helfgott2012}\cite{Helfgott2013}.

\section{Prerequisites}

\subsection{The Discrete Fourier Transform (DFT)}

We will make extensive use of the DFT in what
follows. We adopt the following (un-normalised) definition.

\begin{definition}\label{eq:dft}
Given $N\in\N$ complex values denoted $X_0$ through $X_{N-1}$, the forward
DFT results in $N$ new values $Y_0$ through $Y_{N-1}$
where
\begin{equation*}
Y_m=\sum\limits_{n=0}^{N-1}X_n\e\left(\frac{-nm}{N}\right)
\end{equation*}
and as usual $\e(x):=\exp(2\pi i x)$.

The backward or inverse DFT (iDFT) results from changing the sign in the complex exponential. Performing a forward then backward DFT (or vice versa) multiplies each datum by $N$.
\end{definition}

As written, computing a DFT of length $N$ would appear to have time complexity $\bigoh(N^2)$. The ubiquity of the DFT stems from the existence of $\bigoh(N \log N)$ algorithms, known collectively as the Fast Fourier Transforms
(FFTs). For detailed descriptions of suitable algorithms, we refer the reader
to, for example, \cite{Briggs1995}. However we note that,
significantly for our purposes, this asymptotic complexity can be achieved for arbitrary (even prime) $N$. One such algorithm, and the one we employ, is that due to Bluestein \cite{Bluestein1970}.

Throughout this paper, we will define $F$, the (continuous) Fourier transform of a function $f$ (when it exists), to be
\begin{equation*}
F(x):=\frac{1}{2\pi}\int\limits_{-\infty}^\infty f(t)\exp(-ixt)\D t.
\end{equation*}
Under suitable conditions, the Fourier Inversion Theorem gives us
\begin{equation*}
f(t)=\int\limits_{-\infty}^\infty F(x)\exp(ixt)\D x.
\end{equation*}
To make the transition from the discrete to the continuous, we use the following Theorem.
\begin{theorem}\label{th:dft_pair}
Let $f$ be a function in the Schwartz space with Fourier transform $F$, $N=AB$ with $A,B>0$. Define
\begin{equation*}
\tilde{f}(n):=\sum\limits_{l\in\Z}f\left(\frac{n}{A}+lB\right)
\end{equation*}
 and
\begin{equation*}
\tilde{F}(m):=\sum\limits_{l\in\Z}F\left(\frac{2\pi m}{B}+2\pi lA\right).
\end{equation*}
Then, up to
a constant factor, $\tilde{f}(n)$ and $\tilde{F}(m)$ form a DFT pair of length $N$.
\begin{proof}
By Poisson summation we have
\begin{equation}
\nonumber
\begin{aligned}
\sum\limits_{l\in\Z}f(t+lB)&=\frac{2\pi}{B}\sum\limits_{l\in\Z}F\left(\frac{2\pi l}{B}\right)\e\left(\frac{lt}{B}\right)\\
\tilde{f}(n)&=\frac{2\pi}{B}\sum\limits_{l\in\Z}F\left(\frac{2\pi l}{B}\right)\e\left(\frac{ln}{N}\right).
\end{aligned}
\end{equation}
We now write $l=l^{'}N+m$ to get
\begin{equation}
\nonumber
\begin{aligned}
\tilde{f}(n)&=\frac{2\pi}{B}\sum\limits_{m=0}^{N-1}\sum\limits_{l^{'}\in\Z}F\left(\frac{2\pi(l^{'}N+m)}{B}\right)\e\left(\frac{(l^{'}N+m)n}{N}\right)\\
&=\frac{2\pi}{B}\sum\limits_{m=0}^{N-1}\e\left(\frac{mn}{N}\right)\tilde{F}(m).
\end{aligned}
\end{equation}
This is by definition an iDFT.
\end{proof}
\end{theorem}

The utility of this theorem will be apparent when $f$ and $F$ both decay quickly enough to allow $\tilde{f}(n)$ and $\tilde{F}(m)$ to be approximated by $f\left(\frac{n}{A}\right)$ and $F\left(\frac{m}{B}\right)$ respectively.

\subsection{Interval Arithmetic}

Like Rumely, we chose to manage rounding and truncation errors throughout our
computations using interval arithmetic. We refer the interested reader to the
extensive literature on this subject (perhaps \cite{Moore1966} is a good starting
point) but we summarise our approach below.

Almost all real numbers cannot be represented by a
floating point number of any given precision. Thus, whenever an operation is
carried out on floating point numbers, unless we are very lucky, the answer
will not be exactly representable. We typically attempt to round to the
nearest real number that is exactly representable and thus incur a rounding error. Such errors will
accumulate over time and, to quote Moore
``it is often prohibitively difficult to tell in advance of a
computation how many places must be carried to guarantee results of required
accuracy.'' \cite{Moore1965}.

Instead, we store our intermediate results as two exactly representable
floating point numbers representing an interval that brackets the true
result. The usual mathematical operators and functions are then abstracted to
handle this new data type.

For high precision work (more than the $53$ bits of IEEE double precision
\cite{IEEE754}) we use Revol and Rouillier's MPFI package \cite{Revol2002}. For
computations where double precision will suffice, we use our own
implementation based on the work of Lambov \cite{Lambov2008} for $+$, $-$,
$\times$,$ \div$ and $\rm{sqrt}$. For $\exp$, $\log$, $\sin$, $\cos$ and
$\atan$ we use Muller and de Dinechin's ``Correctly Rounded Mathematical
Library'' \cite{Muller2010}. In both the high precision and double precision
cases, we extend the real interval data type to the complexes in the
obvious (and very probably sub-optimal) way, representing complex values as
rectangles whose corners are exactly representable.

\section{Turing's Method}

Armed with the completed L-function, we have reduced the problem of locating simple zeros of $L_\chi$ on the half line to that of finding sign changes of $\Lambda_\chi$. However, we now need a reference to confirm that all the expected zeros are accounted for. We use a variation on Turing's method from \cite{Turing1953}, extended by Rumely and Trudgian. We start with the following Theorem.
\begin{theorem}[Booker]\label{th:booker}
Let $L(s)$ be an L-function given by an Euler product of degree $r$ and absolutely convergent for $\Re s>1$. Define
\begin{equation*}
\Gamma_\Re(s):=\pi^{-s/2}\Gamma\left(\frac{s}{2}\right),
\end{equation*}
\begin{equation*}
\gamma(s):=\epsilon N^{\frac{1}{2}\left(s-\frac{1}{2}\right)}\prod\limits_{j=1}^r\Gamma_\Re(s+\mu_j),
\end{equation*}
\begin{equation*}
\Lambda(s):=\gamma(s)L(s),
\end{equation*}
where $|\epsilon|=1$, $N\in\Z_{>0}$ and $\Re \mu_j\geq -\frac{1}{2}$ are chosen so that $\Lambda$ satisfies the functional equation
\begin{equation*}
\Lambda(s)=\overline{\Lambda(1-\overline{s})}.
\end{equation*}

Now define
\begin{equation*}
\Phi(t):=\frac{1}{\pi}\left[\arg\epsilon+\frac{\log N}{2}t-\frac{\log\pi}{2}\left(rt+\Im\sum\limits_{j=1}^r\mu_j\right)+\Im\sum\limits_{j=1}^r\log\Gamma\left(\frac{1/2+it+\mu_j}{2}\right)\right]
\end{equation*}
and for $t$ not the ordinate of a zero nor pole of $\Lambda$ define
\begin{equation*}
S(t):=\frac{1}{\pi}\int\limits_\infty^{1/2}\frac{L'}{L}(\sigma+it)\D\sigma.
\end{equation*}
Where $t$ is the ordinate of a zero or pole, set $S(t)=\lim_{\delta\rightarrow 0^+}S(t+\delta)$ (i.e.~ $S$ is upper semi-continuous).
Finally, define
\begin{equation*}
N(t):=\Phi(t)+S(t).
\end{equation*}
Then for $t_1<t_2$, the net number of zeros with imaginary part in $[t_1,t_2)$ counting multiplicity is $N(t_2)-N(t_1)$.
\begin{proof}
See \S $4$ of \cite{Booker2006}.
\end{proof}
\end{theorem}

\begin{theorem}
Given $t_0,h>0$ such that neither $t_0$ nor $t_0+h$ is the imaginary part of a zero of $L_\chi(s)$, let $N_\chi(t_0)$ be the number of zeros, counted with multiplicity, of $L_\chi(s)$ with $|\Im(s)|\leq t_0$ and $\Re(s)\in(0,1)$. Let $\widetilde{N}_{t_0,\chi}(t)$ count the zeros of $L_\chi(s)$ with $\Im(s)\in[t_0,t)$, starting at $0$ at $t_0$ and increasing by $1$ at every zero.

Now for $t$ not the ordinate of a zero of $L_\chi$, define $S_\chi(t)$ by
\begin{equation}
\nonumber
S_\chi(t):=\frac{1}{\pi}\Im \int\limits_\infty^\frac{1}{2}\frac{L_\chi^{'}}{L_\chi}(\sigma+it)\D \sigma
\end{equation}
and take $S_\chi(t)$ to be upper semi-continuous. Then we have
\begin{equation}
\nonumber
\begin{aligned}
N_\chi(t_0)=&\frac{1}{h\pi}\left[2h+\frac{2ht_0+h^2}{2}\log\left(\frac{q}{\pi}\right)+2\int\limits_{t_0}^{t_0+h}\Im\log\Gamma\left(\frac{1/2+a_\chi+it}{2}\right)dt\right.\\
&\left.-\int\limits_{t_0}^{t_0+h}\widetilde{N}_{t_0,\chi}(t)\D t-\int\limits_{t_0}^{t_0+h}\widetilde{N}_{t_0,\overline{\chi}}(t)dt+\int\limits_{t_0}^{t_0+h}S_\chi(t)\D t+\int\limits_{t_0}^{t_0+h}S_{\overline{\chi}}(t)\D t\right].
\end{aligned}
\end{equation}
\begin{proof}
This is Theorem \ref{th:booker} specialised to Dirichlet L-functions. In the terminology of that Theorem, we have $N=q$, $r=1$ and $\mu_1=a_\chi$. We treat conjugate characters in pairs to avoid problems with the arbitrary choice of $\epsilon_\chi$ and to allow for the possibility that $S_\chi(0)$ isn't small. Finally, we integrate both sides from $t_0$ to $t_0+h$.
\end{proof}
\end{theorem}

\begin{theorem} (Rumely).
For $t_0>50$ and $h>0$
\begin{equation}
\nonumber
\left|\int\limits_{t_0}^{t_0+h}S_\chi(t)\D t\right|\leq1.8397+0.1242\log\left(\frac{q(t_0+h)}{2\pi}\right).
\end{equation}
\begin{proof}
Theorem 2 of \cite{Rumely1993}.
\end{proof}
\end{theorem}

Trudgian considered this problem in \cite{Trudgian2011} and in a personal communication, provided revised constants optimised for $qt_0$ in the region of $10^8$. These are $2.17618$ and $0.0679955$ respectively.

\section{An Algorithm for Large $q$}

For ``large'' moduli, we compute the values of $L_\chi(s)$ simultaneously
for all characters of a given modulus by expressing the calculations as a
Discrete Fourier Transform. Specifically, we appeal to the following lemma.

\begin{lemma}\label{lem:dc_dft}
For $q\in\Z\geq 3$ and given $\varphi(q)$ complex values $a(n)$ for $n\in[1,q-1]$ and $(n,q)\neq 0$, we can compute
\begin{equation}%\label{eq:dc_sum}
\nonumber
\sum\limits_{n=1}^{q-1}a(n)\chi(n)
\end{equation}
for the $\varphi(q)$ characters $\chi$ in $\bigoh(\varphi(q)\log(q))$ time and $\bigoh(\varphi(q))$ space.
\begin{proof}
Let $U(R)$ be the group of units of the ring $R$. For $q\in\N$ with the prime decomposition $q=2^\alpha\prod\limits_{i=1}^m p_i^{\alpha_i}$. We consider four cases;
\begin{enumerate}
\item $\alpha=0$ ($q$ is odd) then by the Chinese Remainder Theorem (CRT) we have the constructive, canonical group isomorphism
\begin{equation}
\nonumber
U(\Z/ q\Z)\cong \prod\limits_{i=1}^m U(\Z/ p_i^{\alpha_i}\Z).
\end{equation}
Each of these groups is cyclic so given a primitive root for each $p_i^{\alpha_i}$ we have our construction. Thus this case reduces to performing $\varphi(q)/\varphi(p_i^{\alpha_i})$ length $\varphi(p_i^{\alpha_i})$ DFTs for $i=1\ldots m$.
\item $\alpha=1$ then by the CRT we have the constructive group isomorphism
\begin{equation}
\nonumber
U(\Z/ q\Z)\cong U(\Z/ 2p_1^{\alpha_1}\Z)\prod\limits_{i=2}^m U(\Z/ p_i^{\alpha_i}\Z).
\end{equation}
Each of these groups is cyclic so given a primitive root for $2p_1^{\alpha_1}$ and each $p_i^{\alpha_i}$ ($i>1$) we have our construction. Thus this case reduces to performing $\varphi(q)/\varphi(2p_1^{\alpha_1})$ length $\varphi(2p_1^{\alpha_1})$ DFTs followed by $\varphi(q)/\varphi(p_i^{\alpha_i})$ length $\varphi(p_i^{\alpha_i})$ DFTs for $i=2\ldots m$.
\item $\alpha=2$ then by the CRT we have the constructive, canonical group isomorphism
\begin{equation}
\nonumber
U(\Z/ q\Z)\cong U(\Z/ 4\Z)\prod\limits_{i=1}^m U(\Z/ p_i^{\alpha_i}\Z).
\end{equation}
Each of these groups is cyclic so given a primitive root for each $p_i^{\alpha_i}$ ($i>1$) we have our construction. Thus this case reduces to performing $\varphi(q)/2$ length $2$ DFTs followed by $\varphi(q)/\varphi(p_i^{\alpha_i})$ length $\varphi(p_i^{\alpha_i})$ DFTs for $i=1\ldots m$.
\item $\alpha>2$ then by the CRT we have the constructive, canonical group isomorphism
\begin{equation}
\nonumber
U(\Z/ q\Z)\cong U(\Z/ 2^\alpha\Z)\prod\limits_{i=1}^m U(\Z/ p_i^{\alpha_i}\Z).
\end{equation}
Now $U(\Z/ 2^\alpha\Z)$ is the product of a cyclic group of order $2$ and a cyclic group of order $2^{\alpha-2}$ with pseudo primitive roots $-1$ and $5$ respectively. The remaining groups (if there are any) are cyclic so given a primitive root for each $p_i^{\alpha_i}$ ($i>1$) we have our construction. Thus this case reduces to performing $\varphi(q)/2$ length $2$ DFTs, $\varphi(q)/2^{\alpha-2}$ length $2^{\alpha-2}$ DFTs followed by $\varphi(q)/\varphi(p_i^{\alpha_i})$ length $\varphi(p_i^{\alpha_i})$ DFTs for $i=1\ldots m$.
\end{enumerate}
In each case, given the ability to perform an arbitrary length $n$ DFT in time $\bigoh(n\log n)$, we have the claimed overall complexity.
\end{proof}
\end{lemma}

We seek to apply Lemma \ref{lem:dc_dft} by way of the Hurwitz zeta function, defined for $\Re s>1$ and $\alpha\in(0,1]$ by
\begin{equation*}
\zeta(s,\alpha):=\sum\limits_{n=0}^\infty (n+\alpha)^{-s}.
\end{equation*}

This function has analytic continuation to $\C$ with the exception of a simple
pole at $s=1$ and except at this pole it can be used to express any Dirichlet
L-function of modulus $q$ via
\begin{equation*}
L_\chi(s)=q^{-s}\sum\limits_{a=1}^q\chi(a)\zeta\left(s,\frac{a}{q}\right)
\end{equation*}
(see $\S$ 12 of \cite{Apostol1976}).

Thus, for a given $q$ and $s$, if we can supply the $\varphi(q)$ values of
$\zeta\left(s,\frac{a}{q}\right)$ for $a\in[1,q-1]$ with $(a,q)=1$, we can apply
Lemma \ref{lem:dc_dft} to compute each $L_\chi(s)$ in, on average, time $\log q$.

\subsection{Computing $\zeta(1/2+it,a/q)$}

For a given $t\in\R_{\geq 0}$ and $q\geq 3$, we need to be able to compute $\zeta(1/2+it,a/q)$ for $a\in[1,q)$ with
$(a,q)=1$. We proceed by computing, for each $t$, a lattice of $D$ rows and $N$ columns were the entry in the $r$'th row and $c$'th column are $\zeta(1/2+it+c,r/D)$ ($r$ running $1\ldots D$ and $c$ $0\ldots N$). We chose $N=15$ and $D=2,048$ to achieve the necessary precision.

We use the following lemma both to initially compute the lattice (once, in high precision using MPFI) and to compute the required values for the DFT from that lattice (many times, using double precision intervals).
 
\begin{lemma}\label{lem:hur_tay}
For $s\not\in\Z_{\leq 0}$, $\alpha\in(0,1]$ and $|\delta|<\alpha$
\begin{equation}
%\label{eq:hurtaylor}
\nonumber
\zeta(s,\alpha+\delta)=\sum\limits_{k=0}^\infty \frac{(-\delta)^k\zeta(s+k,\alpha)\prod_{j=0}^{k-1}(s+j)}{k!}.
\end{equation}
\begin{proof}
Starting with $\Re s>1$ and differentiating term by term we have
\begin{equation}
\nonumber
\zeta^{(k)}(s,\alpha)=\sum\limits_{n=0}^{\infty}(-1)^ks(s+1)(s+2)\ldots(s+k-1)(n+\alpha)^{-s-k}
\end{equation}
and the result follows for $\Re s>1$ by Taylor's Theorem. The Taylor expansion
also gives us the analytic continuation to $\C\setminus\Z_{\leq0}$.
\end{proof}
\end{lemma}
In practice, it is better to work with
\begin{equation}
\nonumber
\zeta_{M}(s,\alpha)=\zeta(s,\alpha)-\sum\limits_{n=0}^M (n+\alpha)^{-s}
\end{equation}
for some $M\in\N$ and to recover
$\zeta(s,\alpha)$ by adding back the missing terms.
% To be able to rigorously bound the error in truncating the series definition and the Taylor approximation, we use the following lemmas.

% \begin{lemma}
% \label{lem:hur_sumbound}
% For $\alpha\in(0,1]$, $\Re(s)>1$ and $M\geq2$
% \begin{equation}
% \nonumber
% \left|\zeta_{M}(s,\alpha)\right|\leq\frac{(M+\alpha-1)^{1-\Re(s)}}{\Re(s)-1}.
% \end{equation}
% \begin{proof}
% Integral test.
% \end{proof}
% \end{lemma}

% \begin{lemma}
% \label{lem:hur_taybound}
% If we use the first $N$ terms of the Taylor approximation to $\zeta_{M}(s,\alpha+\delta)$, then the absolute error is bounded by
% \begin{equation}
% \nonumber
% \frac{(N+1)\left|s(s+1)\ldots (s+N-1)\zeta_{M}(s+N,\alpha)\right|\delta^N}{N!(N+1-(|s|+N)\delta)}
% \end{equation}
% and the approximation is valid for
% \begin{equation}
% \nonumber
% \frac{(|s|+N)\delta}{N+1}<1.
% \end{equation}
% \begin{proof}
% The first term dropped is
% \begin{equation}
% \nonumber
% \frac{s(s+1)\ldots (s+N-1)\zeta_{M}(s+N,\alpha)\delta^N}{N!}
% \end{equation}
% and the result follows by considering the geometric sequence with this first term and with common ratio
% \begin{equation}
% \nonumber
% \frac{(|s|+N)\delta}{N+1}.
% \end{equation}
% \end{proof}
% \end{lemma}

\section{An Algorithm for Small $q$}

The algorithm described above starts to become unwieldy as $t$, the height up the critical line, increases. Each new $t$ requires its own pre-computed lattice of $\zeta(1/2+it+c,r/D)$ and the cost of producing this lattice is amortised over less and less $q$. In \cite{Booker2006}, Booker describes a rigorous algorithm for computing L-functions. What follows is that algorithm specialised to Dirichlet L-functions.

For $\eta\in(-1,1)$ and even primitive characters $\chi$ define
\begin{equation}
\begin{aligned}
\nonumber
F_e(t,\chi):=&\epsilon_\chi q^{\frac{it}{2}} \pi^{-\frac{1/2+it}{2}}
\Gamma\left(\frac{1/2+it}{2}\right) \exp\left(\frac{\pi\eta
    t}{4}\right)L_\chi\left(\frac{1}{2}+it\right)\textrm{ and}\\
\hat{F}_e(x,\chi):=&\frac{1}{2\pi}\int\limits_{-\infty}^{\infty}F_e(t,\chi)e^{-ixt}\D t.
\end{aligned}
\end{equation}
For odd primitive characters $\chi$ define 
\begin{equation}
\begin{aligned}
\nonumber
F_o(t,\chi):=&\epsilon_\chi q^{\frac{it}{2}} \pi^{-\frac{3/2+it}{2}}
\Gamma\left(\frac{3/2+it}{2}\right) \exp\left(\frac{\pi\eta
    t}{4}\right)L_\chi\left(\frac{1}{2}+it\right)\textrm{ and}\\
\hat{F}_o(x,\chi):=&\frac{1}{2\pi}\int\limits_{-\infty}^{\infty}F_o(t,\chi)e^{-ixt}\D t.
\end{aligned}
\end{equation}
We chose the parameter $\eta$ to control the decay of the gamma factor as $t$
increases.

We now choose $A,B>0$ with $AB\in2^{\N}$ and define
\begin{equation}%\label{eq:fhattwiddle}
\nonumber
\begin{aligned}
\widetilde{\hat{F}}_e(n,\chi)&:=\sum\limits_{k\in\Z}\hat{F}_e\left(\frac{2\pi
    n}{B}+2\pi kA,\chi\right)\\
\textrm{and}\\
\widetilde{\hat{F}}_o(n,\chi)&:=\sum\limits_{k\in\Z}\hat{F}_o\left(\frac{2\pi
    n}{B}+2\pi kA,\chi\right).
\end{aligned}
\end{equation}
Similarly, define
\begin{equation}%\label{eq:ftwiddle}
\nonumber
\begin{aligned}
\widetilde{F}_e(m,\chi)&:=\sum\limits_{k\in\Z}F_e\left(\frac{m}{A}+kB,\chi\right)\\
\textrm{and}\\
\widetilde{F}_o(m,\chi)&:=\sum\limits_{k\in\Z}F_o\left(\frac{m}{A}+kB,\chi\right).
\end{aligned}
\end{equation}
In outline, the method is
\begin{enumerate}
\item Compute $\hat{F}_e\left(\frac{2\pi n}{B}\right)$ or $\hat{F}_o\left(\frac{2\pi n}{B}\right)$ for $n=0\ldots N-1$.
\item Use these values as an approximation to $\widetilde{\hat{F}}_e(n,\chi)$ or $\widetilde{\hat{F}}_o(n,\chi)$ respectively.
\item Appealing to Theorem \ref{th:dft_pair}, perform a DFT to yield
  $\widetilde{F}_e(m,\chi)$ or $\widetilde{F}_o(m,\chi)$ respectively.
\item Use $\widetilde{F}_e(m,\chi)$ or $\widetilde{F}_o(m,\chi)$ as an approximation to
  $F_e\left(\frac{m}{A},\chi\right)$ or $F_o\left(\frac{m}{A},\chi\right)$ respectively.
\end{enumerate}

%  Continuing in the notation of \cite{Booker2006}, we have $r=1$, $m=0$, $P(s)=1$, $C=1$, $\alpha=0$ and $c'=0$ with the balance of the parameters as per Table \ref{tab:pars}.

% \begin{table}
% \caption{Parameters}
% \label{tab:pars}\centering
% \begin {tabular}{ccc}
% \hline
%  & Even $\chi$ & Odd $\chi$ \\\hline
% $\mu_1$ & $0$ & $1$ \\
% $\nu_1$ & $0$ & $\frac{1}{2}$ \\
% $\mu$ & $\frac{1}{2}$ & $\frac{3}{2}$ \\
% $c$ & $1$ & $2$ \\\hline           
% \end{tabular}
% \end{table}

We now make the above outline rigorous.

\subsection{Computing $\hat{F}_e(t)$ and $\hat{F}_o(t)$}

\begin{lemma}\label{lem:fehat}
Let $x\in\R$, $\eta\in(-1,1)$ and $u(x):=\frac{\pi\eta i}{4}+x$. Then we have
\begin{equation}
\nonumber
\hat{F}_e(x,\chi)=\frac{2\epsilon_\chi\exp\left(\frac{u(x)}{2}\right)}{q^\frac{1}{4}}\sum\limits_{n=1}^\infty \chi(n)\exp\left(-\frac{\pi n^2\exp(2u(x))}{q}\right).
\end{equation}
\begin{proof}
Writing $s=1/2+it$ we get
\begin{equation}
\nonumber
\begin{aligned}
\hat{F}_e(x,\chi)=&\frac{\epsilon_\chi}{2\pi i}\int\limits_{\Re(s)=\frac{1}{2}}q^\frac{s-1/2}{2}\pi^{-\frac{s}{2}}\Gamma\left(\frac{s}{2}\right)\exp\left(\frac{-\pi \eta i(s-1/2)}{4}\right)\exp(-x(s-1/2))L_\chi(s)\D s\\
=&\frac{\epsilon_\chi}{2\pi i}\int\limits_{\Re(s)=2}q^\frac{s-1/2}{2}\pi^{-\frac{s}{2}}\Gamma\left(\frac{s}{2}\right)\exp\left(\frac{-\pi \eta i(s-1/2)}{4}\right)\exp(-x(s-1/2))L_\chi(s)\D s\\
=&\frac{\epsilon_\chi}{q^\frac{1}{4}}\frac{1}{2\pi i}\int\limits_{\Re(s)=2} \left(\frac{q}{\pi}\right)^\frac{s}{2}\Gamma\left(\frac{s}{2}\right)\exp\left(-\left(\frac{\pi \eta i+4x}{4}\right)(s-1/2)\right)\sum\limits_{n=1}^\infty \chi(n)n^{-s}\D s\\
=&\frac{\epsilon_\chi\exp(u(x)/2)}{q^\frac{1}{4}}\sum\limits_{n=1}^\infty \chi(n)\frac{1}{2\pi i}\int\limits_{\Re(s)=2} \left(\frac{\pi n^2}{q}\right)^{-\frac{s}{2}}\Gamma\left(\frac{s}{2}\right)\exp(2u(x))^{-s/2}\D s\\
=&\frac{2\epsilon_\chi\exp\left(\frac{u(x)}{2}\right)}{q^\frac{1}{4}}\sum\limits_{n=1}^\infty \chi(n)\exp\left(-\frac{\pi n^2\exp(2u(x))}{q}\right).
\end{aligned}
\end{equation}
\end{proof}
\end{lemma}

We can rigorously bound the error in truncating the sum either by reference to Lemma $5.4$ of \cite{Booker2006} or by majorising the missing terms with the obvious geometric series.

\begin{lemma}\label{lem:fohat}
Let $x$, $\eta$ and $u(x)$ be as defined in Lemma \ref{lem:fehat}. Then we have 
\begin{displaymath}
\nonumber
\hat{F}_o(x,\chi)=\frac{2\epsilon_\chi\exp\left(\frac{3u(x)}{2}\right)}{q^\frac{3}{4}}\sum\limits_{n=1}^\infty n\chi(n)\exp\left(-\frac{\pi n^2\exp(2u(x))}{q}\right).
\end{displaymath}
\begin{proof}
The proof follows the same lines as Lemma \ref{lem:fehat}.
\end{proof}
\end{lemma}

\subsection{Approximating $\widetilde{\hat{F}}_e$ and $\widetilde{\hat{F}}_o$ with $\hat{F}_e$ and $\hat{F}_o$}

We intend to chose our parameters to allow us to use $\hat{F}_e$ and $\hat{F}_o$ as approximations to $\widetilde{\hat{F}}_e$ and $\widetilde{\hat{F}}_o$ respectively. We therefore need to bound the error introduced and we start with two lemmas.
\begin{lemma}
\label{lem:lfuncbound}
For $t\in\R$ we have
\begin{equation}
\nonumber
\left|L_\chi\left(\frac{1}{2}+it\right)\right|\leq\zeta\left(\frac{9}{8}\right)\left(\frac{q}{2\pi}\right)^{5/16}\left(\frac{3}{2}+|t|\right)^{5/16}.
\end{equation}
\begin{proof}

We evaluate Rademacher's bound \cite{Rademacher1959}
\begin{equation}
\nonumber
\left|L_\chi(s)\right|\leq\zeta(1+\nu)\left(\frac{q|1+s|}{2\pi}\right)^{\frac{1+\nu-\Re(s)}{2}}
\end{equation}
 with $\nu=1/8$ and $s=1/2+it$.
\end{proof}
\end{lemma}

\begin{lemma}{(Booker)}\label{lem:bookerfhat}
Let $\eta\in(-1,1)$, $\delta=\frac{\pi}{2}(1-|\eta|)$ and $X(x)=\pi\delta\e^{-\delta}x>1$. Then
\begin{equation*}
\left|\sum\limits_{k=0}^\infty \hat{F}_e(x+2\pi kA)\right|\leq\frac{4\exp\left(\frac{x}{2}-X(x)\right)\left(1+\frac{1}{2X(x)}\right)}{\delta^{\frac{1}{2}}q^{\frac{1}{4}}\left(1-\e^{-\pi A}\right)}
\end{equation*}
and
\begin{equation*}
\left|\sum\limits_{k=0}^\infty \hat{F}_o(x+2\pi kA)\right|\leq\frac{4\exp\left(\frac{3x}{2}-X(x)\right)\left(1+\frac{1}{2X(x)}\right)^{\frac{3}{2}}}{\delta^{\frac{1}{2}}q^{\frac{3}{4}}\left(1-\e^{-\pi A}\right)}.
\end{equation*}

\begin{proof}
This is Lemma $5.6$ of \cite{Booker2006} specialised to Dirichlet L-functions.
\end{proof}
\end{lemma}

We can now proceed to the necessary bounds.

\begin{lemma}
Let $A\geq\frac{1}{2\pi}$, $B>0$, $w_1=\frac{2\pi n}{B}+2\pi A$, $w_2=-\frac{2\pi n}{B}+2\pi A$, with $X(x)$ and $\delta$ as defined in Lemma \ref{lem:bookerfhat} and $X(w_1),X(w_2)>1$. Then
\begin{displaymath}
\left|\widetilde{\hat{F}}_e(n,\chi)-\hat{F}_e\left(\frac{2\pi n}{B},\chi\right)\right|\leq
\end{displaymath}
\begin{displaymath}
\frac{4\left(\exp\left(\frac{w_1}{2}-X(w_1)\right)\left(1+\frac{1}{2X(w_1)}\right)\right.+\left.\exp\left(\frac{w_2}{2}-X(w_2)\right)\left(1+\frac{1}{2X(w_2)}\right)\right)}{q^{1/4}\delta^{1/2}(1-e^{-\pi A})}
\end{displaymath}
and
\begin{displaymath}
\left|\widetilde{\hat{F}}_o(n,\chi)-\hat{F}_o\left(\frac{2\pi n}{B},\chi\right)\right|\leq
\end{displaymath}
\begin{displaymath}
\frac{4\left(\exp\left(\frac{3w_1}{2}-X(w_1)\right)\left(1+\frac{1}{2X(w_1)}\right)^\frac{3}{2}\right.+\left.\exp\left(\frac{3w_2}{2}-X(w_2)\right)\left(1+\frac{1}{2X(w_2)}\right)^\frac{3}{2}\right)}{q^{3/4}\delta^{1/2}(1-e^{-\pi A})}.
\end{displaymath}

\begin{proof}
We apply Lemma \ref{lem:bookerfhat} with $x=\frac{2\pi n}{B}\pm2\pi A$.
\end{proof}
\end{lemma}

\begin{lemma}
Given $t\in\R$ and $B>0$, we define
\begin{displaymath}
E_e(t):=\zeta\left(\frac{9}{8}\right)\pi^{-\frac{1}{4}}\left|\Gamma\left(\frac{1}{4}+\frac{it}{2}\right)\right|e^{\frac{\pi}{4}\eta t}\left(\frac{q}{2\pi}\left|\frac{3}{2}+t\right|\right)^{\frac{5}{16}},
\end{displaymath}
\begin{displaymath}
\beta_e(t):=\frac{\pi}{4}-\frac{1}{2}\arctan\left(\frac{1}{2|t|}\right)-\frac{4}{\pi^2|t^2-\frac{1}{4}|},
\end{displaymath}
\begin{displaymath}
E_o(t):=\zeta\left(\frac{9}{8}\right)\pi^{-\frac{3}{4}}\left|\Gamma\left(\frac{3}{4}+\frac{it}{2}\right)\right|e^{\frac{\pi}{4}\eta t}\left(\frac{q}{2\pi}\left|\frac{3}{2}+t\right|\right)^{\frac{5}{16}}
\end{displaymath}
and
\begin{displaymath}
\beta_o(t):=\frac{\pi}{4}-\frac{3}{2}\arctan\left(\frac{1}{2|t|}\right)-\frac{4}{\pi^2|t^2-\frac{9}{4}|}.
\end{displaymath}
Then for $\beta_{e,o}\left(\frac{m}{A}+B\right)>\frac{\pi}{4}\eta$ and $\beta_{e,o}\left(\frac{m}{A}-B\right)>-\frac{\pi}{4}\eta$ we have
\begin{displaymath}
\left|\widetilde{F}_e(m,\chi)-F_e\left(\frac{m}{A},\chi\right)\right|\leq
\end{displaymath}
\begin{displaymath}
\frac{E_e\left(\frac{m}{A}+B\right)}{1-\exp(-B(\beta_e(m/A+B)-\frac{\pi}{4}\eta))}+\frac{E_e\left(\frac{m}{A}-B\right)}{1-\exp(-B(\beta_e(m/A-B)+\frac{\pi}{4}\eta))}
\end{displaymath}
and
\begin{displaymath}
\left|\widetilde{F}_o(m,\chi)-F_o\left(\frac{m}{A},\chi\right)\right|\leq
\end{displaymath}
\begin{displaymath}
\frac{E_o\left(\frac{m}{A}+B\right)}{1-\exp(-B(\beta_o(m/A+B)-\frac{\pi}{4}\eta))}+\frac{E_o\left(\frac{m}{A}-B\right)}{1-\exp(-B(\beta_o(m/A-B)+\frac{\pi}{4}\eta))}.
\end{displaymath}
\begin{proof}
We apply Lemma 5.7 (i) of \cite{Booker2006} with $t=\frac{m}{A}+B$ and 5.7 (ii) with $t=\frac{m}{A}-B$, replacing the bound for $L_\chi(s)$ with our Lemma \ref{lem:lfuncbound}.
\end{proof}
\end{lemma}
We note here that the condition on $\beta_{e,o}(t)$ will fail when $t$ is small, i.e.~when $\frac{m}{A}\approx B$. However, this only happens for $m$ approaching $AB$, by which point the loss of precision through other factors has rendered these values useless for computational purposes anyway.

\section{Rigorous Up-sampling}

The output from both algorithms is a lattice of values of $\Lambda_\chi(t)$. The sample rate used ($5/64$) will be insufficient to resolve all the zeros, so we employ a rigorous up-sampling technique based on theorems of Whittaker-Shannon and Weiss.

\begin{theorem}\label{th:Whittaker}
(Whittaker-Shannon Sampling Theorem) Let $f(t)$ be a continuous, real valued function with Fourier Transform $F(x)$ such that $F(x)=0$ for $|x|>2\pi B>0$ (i.e.~$f(t)$ is band-limited with bandwidth $2\pi B$). Also, define
\begin{equation}
\nonumber
\sinc(x):=\frac{\sin(x)}{x}.
\end{equation}
Then
\begin{equation}
\nonumber
f(t)=\sum\limits_{n\in\Z}f\left(\frac{n}{2B}\right)\sinc\left(2B\pi\left(\frac{n}{2B}-t\right)\right),
\end{equation}
when this sum converges.
\begin{proof}
See \cite{Walker1991}.
\end{proof}
\end{theorem}

To apply Theorem \ref{th:Whittaker} rigorously, we need to examine two sources of error
\begin{itemize}
\item the error introduced by truncating the sum and
\item the error introduced if the function is only approximately band-limited.
\end{itemize}

The former will be dealt with on a case by case basis. The latter, referred to as aliasing in signal processing circles, is the subject of a theorem due to Weiss.
\begin{theorem}\label{th:Weiss}
(Weiss) Let $f(t)$ be a real valued function with Fourier Transform $F(x)$ such that
\begin{enumerate}
	\item $\int\limits_{-\infty}^\infty|F(x)|dx < \infty$
	\item $F(x)$ is of bounded variation on $\R$
	\item when $F$ has a jump discontinuity at $x$, then
          $F(x)=\lim\limits_{\epsilon\rightarrow 0^+} \frac{F(x-\epsilon)+F(x+\epsilon)}{2}$.
\end{enumerate}
Then
\begin{equation}%\label{eq:weiss}
\nonumber
\left|f(t)-\sum\limits_{n\in\Z}f\left(\frac{n}{2B}\right)\sinc\left(2B\pi\left(t-\frac{n}{2B}\right)\right)\right|\leq 4\int\limits_{2\pi B}^{\infty}\left|F(x)\right|dx.
\end{equation}
\begin{proof}
See for example \cite{Brownjr.1967}.
\end{proof}
\end{theorem}

For $t_0\in\R$ and $h>0$ define $W:\R\rightarrow\R$ by
\begin{equation}
\nonumber
W(t,\chi):=\Lambda_\chi(t)\exp\left(\frac{-(t-t_0)^2}{2h^2}\right)
\end{equation}

so $W(t_0,\chi)=\Lambda_\chi(t)$.

We aim to estimate $W(t_0,\chi)$ from our samples using Theorems \ref{th:Whittaker} (Whittaker-Shannon) and \ref{th:Weiss}. The following lemmas provide the necessary rigorous bounds.

\begin{lemma}
\label{lem:gammabound}
For $a_\chi\in\{0,1\}$
\begin{equation}
\nonumber
\begin{aligned}
&\left|\Gamma\left(\frac{\frac{1}{2}+it+a_\chi}{2}\right)\right|e^\frac{\pi t}{4}\\
&\;\;\;\;\;\;\;\leq\max\left(2^{1/4}\sqrt{\pi}\left(\frac{3}{2}+\max(t,0)\right)^{\frac{1}{4}}\exp\left(\frac{1}{6}\right),\sqrt{2\pi}\exp\left(\frac{\pi}{8}+\frac{1}{4}\right)\right).
\end{aligned}
\end{equation}
\begin{proof}
We use Stirling's approximation separately for $a_\chi=0$ and $a_\chi=1$.
\end{proof}
\end{lemma}

\begin{lemma}
\label{lem:LIbound}
Define $I_\chi$ by 
\begin{equation}
\nonumber
I_\chi:=\frac{2}{\pi}\int\limits_{2\pi B}^\infty\left|\int\limits_{-\infty}^\infty W(t,\chi)\exp(-ixt)\D t\right|\D x.
\end{equation}

Then, writing $M$ in place of $\frac{5}{2}-a_\chi$ we have 
\begin{equation}
\nonumber
I_\chi\leq\frac{2\left(\frac{q}{\pi}\right)^\frac{M}{2}\zeta(M+1/2)\exp\left(\frac{M^2}{2h^2}-2\pi BM\right)P(t_0,h)}{\pi M}
\end{equation}
where
\begin{equation}
\nonumber
P(t_0,h)=\int\limits_{-\infty}^{\infty} \left|\Gamma\left(\frac{3+it}{2}\right)\right|\exp\left(\frac{\pi t}{4}\right)\exp\left(-\frac{(t-t_0)^2}{2h^2}\right)\D t.
\end{equation}
\begin{proof}
Writing $s=1/2+it$ we get
\begin{equation}
\nonumber
\begin{aligned}
I_\chi\leq&\frac{2}{\pi}\int\limits_{2\pi B}^\infty\left|\int\limits_{\Re(s)=1/2} \left(\frac{q}{\pi}\right)^{\frac{s-1/2}{2}}\Gamma\left(\frac{s+a_\chi}{2}\right)\exp\left(\frac{\pi i(1/2-s)}{4}\right)L_\chi(s)\right.\\
&\;\;\;\;\;\;\;\;\;\;\;\;\;\;\;\;\;\;\;\;\;\;\;\left.\exp((1/2-s)x)\exp\left(\frac{-(i(1/2-s)-t_0)^2}{2h^2}\right)\D s\right|\D x.
\end{aligned}
\end{equation}
We now shift the contour of integration to the right so that $\Re(s)=\sigma=3-a_\chi$ and write $s=M+1/2+it$ to get
\begin{equation}
\nonumber
\begin{aligned}
I_\chi\leq\frac{2}{\pi}\int\limits_{2\pi B}^\infty\int\limits_{-\infty}^\infty&\left| \left(\frac{q}{\pi}\right)^{\frac{M}{2}}\Gamma\left(\frac{3+it}{2}\right)\exp\left(\frac{\pi t}{4}\right)\zeta(M+1/2)\right.\\
&\;\;\;\left.\exp(-Mx)\exp\left(\frac{M^2-(t-t_0)^2}{2h^2}\right)\right|\D t\;\D x.
\end{aligned}
\end{equation}
Integrating with respect to $t$ gives us
\begin{equation}
\nonumber
\begin{aligned}
I_\chi\leq\frac{2}{\pi}&\left(\frac{q}{\pi}\right)^{\frac{M}{2}}\zeta(M+1/2)\exp\left(\frac{M^2}{2h^2}\right)P(t_0,h)\int\limits_{2\pi B}^\infty\exp(-Mx)\D x
\end{aligned}
\end{equation}
and the result follows after integrating with respect to $x$.
\end{proof}
\end{lemma}

\begin{lemma}
\label{lem:Lpmbound}
Let $t_0\geq 0$. Then
\begin{equation}
\nonumber
P(t_0,h)\leq h\pi\left(t_0+\frac{h}{\sqrt{2\pi}}+1+\frac{1}{2\sqrt{2}}\right).
\end{equation}
\begin{proof}
We have
\begin{equation}
\nonumber
\begin{aligned}
P(t_0,h)\leq&\int\limits_0^\infty\left|\Gamma\left(\frac{3+it}{2}\right)\right|\exp\left(\frac{\pi t}{4}\right)\exp\left(\frac{-(t-t_0)^2}{2h^2}\right)\D t\\
+&\int\limits_{-\infty}^0\left|\Gamma\left(\frac{3+it}{2}\right)\right|\exp\left(\frac{\pi t}{4}\right)\exp\left(\frac{-(t-t_0)^2}{2h^2}\right)\D t\\
\leq&\int\limits_0^\infty\frac{(1+t)}{2}\left|\Gamma\left(\frac{1+it}{2}\right)\right|\exp\left(\frac{\pi t}{4}\right)\exp\left(\frac{-(t-t_0)^2}{2h^2}\right)\D t\\
&\;\;\;\;\;\;\;+\Gamma\left(\frac{3}{2}\right)\frac{h\sqrt{2\pi}}{2}\left(1-\erf\left(\frac{\sqrt{2}t_0}{2}\right)\right)\\
\leq&\int\limits_0^\infty\frac{(1+t)}{2}\sqrt{\frac{\pi}{\cosh(\pi t/2)}}\exp\left(\frac{\pi t}{4}\right)\exp\left(\frac{-(t-t_0)^2}{2h^2}\right)dt\\
&\;\;\;\;\;\;\;+\Gamma\left(\frac{3}{2}\right)\frac{h\sqrt{2\pi}}{2}\\
\leq&\int\limits_0^\infty\frac{(1+t)}{2}\sqrt{2\pi}\exp\left(\frac{-(t-t_0)^2}{2h^2}\right)\D t+\frac{h\pi\sqrt{2}}{4}\\
\leq&h\pi\left(\frac{h}{\sqrt{2\pi}}+t_0+1\right)+\frac{h\pi\sqrt{2}}{4}.
\end{aligned}
\end{equation}
\end{proof}
\end{lemma}

\begin{lemma}
\label{lem:gaussbound}
Let $h,B>0$, $t_0=\frac{n_0}{2B}$ for some $n_0\in\N$ and $N\in\N$. Now define
\begin{equation}
\nonumber
\begin{aligned}
G(n):=\frac{\left(\frac{3}{2}+t_0+\frac{N+n}{2B}\right)^{9/16}\exp\left(\frac{-(N+n)^2}{8B^2h^2}\right)}{\pi (N+n)}.
\end{aligned}
\end{equation}
Then
\begin{equation}
\nonumber
\begin{aligned}
\sum\limits_{n\geq2Bt_0+N}&\left(\frac{3}{2}+\frac{n}{2B}\right)^{9/16}\exp\left(\frac{-\left(\frac{n}{2B}-t_0\right)^2}{2h^2}\right)\sinc\left(2B\pi\left(\frac{n}{2B}-t_0\right)\right)\\
&\leq\frac{G(0)}{1-G(1)/G(0)}.
\end{aligned}
\end{equation}
\begin{proof}
$G(n)$ is at least as large as the corresponding term in the sum and the ratio $G(n+1)/G(n)$ is a decreasing function of $n$ so the result follows as the sum of a geometric series.
\end{proof}
\end{lemma}

We can now combine Lemmas \ref{lem:lfuncbound}, \ref{lem:gammabound} and \ref{lem:gaussbound}.
\begin{lemma}
Define
\begin{equation}
\nonumber
E:=\sum\limits_{|n|\geq N}W\left(\frac{n}{2B}\right)\sinc\left(2B\pi\left(\frac{n}{2B}-t_0\right)\right).
\end{equation}
Then for large enough $t_0$ we have
\begin{equation}
\nonumber
|E|\leq\sqrt{\pi}\zeta\left(\frac{9}{8}\right)\exp(1/6)2^{5/4}\left(\frac{q}{2\pi}\right)^{5/16}\frac{G(0)}{1-G(1)/G(0)}.
\end{equation}
\end{lemma}

\section{Results}

Both algorithms parallelise trivially and we ran both algorithms on various clusters in the UK and France. We ensured that every modulus was checked at some point on a system benefiting from ECC memory using the small $q$ algorithm up to $q=10,000$ or so, and the large $q$ algorithm beyond that. We moved up the critical line in steps of $\frac{5}{64}$ representing a sampling rate of about $5$ times the expected zero density. We then routinely up-sampled by a factor of $8$ and then if necessary by $32$, $128$ and ultimately $512$. At this point, about $0.0003\%$ of the L-functions remained due to one or more of the following issues:-
\begin{itemize}
\item The sign of $\Lambda_\chi(1/2)$ could not be determined. This was resolved using a double precision interval implementation of Euler-MacLaurin.
\item The sign of $\Lambda_\chi$ was positive, became indeterminate and then became positive again (or negative, indeterminate, negative). Since a failure to cross the $x$ axis here would, on its own, be enough to refute GRH, we fully expected to find that the indeterminate region was actually hiding a pair of zeros. In every case, using an interval arithmetic version of Euler MacLaurin (first at double precision, but occasionally resorting to MPFI  at $100$ bits) located the expected sign changes.
\item The sign of $\Lambda_\chi$ was positive, indeterminate and then negative (or vice versa). Rather than hiding a single sign change, closer inspection revealed three sign changes in the indeterminate region.
\item Occasionally, the estimate for the number of zeros to locate computed via Turing's method did not bracket an integer. This was caused by zeros being missed in the region used to compute the Turing estimate itself and these were resolved by shifting the region or locating the missing zeros using high precision.
\end{itemize}

In all, the computation consumed approximately $400,000$ core hours.\footnote{The computing resources used were Intel/AMD based and equipped with the SSE2 instruction set. Except for small $q$, where the lengths of the FFTs involved became the limiting factor, we were able to exploit all of the cores available to us on multi-core systems.} We checked all the $29,565,923,837$ Dirichlet L-functions with primitive modulus $q\leq 400,000$, isolating approximately $3.8\cdot 10^{13}$ zeros (not counting those used in Turing's method). Specifically, we have;

\begin{theorem}
GRH holds for Dirichlet L-functions of primitive character modulus $q\leq 400,000$ and to height $t_0=\textrm{max}\left(\frac{10^8}{q},\frac{7.5\cdot 10^7}{q}+200\right)$ for even $q$ and to height $t_0=\textrm{max}\left(\frac{10^8}{q},\frac{3.75\cdot 10^7}{q}+200\right)$ for odd $q$.
\end{theorem}

In addition, we explored the central point of the $739,151,526,102$ primitive characters with $q\leq 2,000,000$ using the large $q$ algorithm. In $438,152$ cases, the computation returned a value for the completed L-function as a double precision interval that straddled zero. Recomputing these points, again using double precision intervals but this time via Euler-MacLaurin, resolved all but $20$ and these were in turn eliminated using Euler-MacLaurin implemented in MPFI at $100$ bits of precision. We can therefore state; 

\begin{theorem}
For every Dirichlet L-function of primitive character modulus $q\leq 2,000,000$, we have
\begin{equation*}
L_\chi(1/2)\neq 0.
\end{equation*}
\end{theorem}

\bibliographystyle{amsplain}
\bibliography{davebib5}{}

\providecommand{\bysame}{\leavevmode\hbox to3em{\hrulefill}\thinspace}
\providecommand{\MR}{\relax\ifhmode\unskip\space\fi MR }
\providecommand{\MRhref}[2]{%
  \href{http://www.ams.org/mathscinet-getitem?mr=#1}{#2}
}
\providecommand{\href}[2]{#2}
\begin{thebibliography}{10}

\bibitem{Apostol1976}
Tom~M. Apostol, \emph{{Introduction to Analytic Number Theory}}, Undergraduate
  Texts in Mathematics, Springer, 1976.

\bibitem{Bluestein1970}
L.~Bluestein, \emph{{A linear filtering approach to the computation of discrete
  Fourier transform}}, IEEE Transactions on Audio and Electroacoustics
  \textbf{18} (1970), no.~4, 451--455.

\bibitem{Booker2006}
Andrew~R. Booker, \emph{{Artin's conjecture, Turing's method and the Riemann
  hypothesis}}, Experiment. Math. \textbf{15} (2006), no.~4, 385--407.

\bibitem{Briggs1995}
William~L. Briggs and Van~Emden Henson, \emph{{The DFT: An Owners Manual for
  the Discrete Fourier Transform}}, SIAM, 1995.

\bibitem{Brownjr.1967}
J.L. Brown~Jr., \emph{{On the Error in Reconstructing a Non-Bandlimited
  Function by Means of the Bandpass Sampling Theorem}}, J. Math. Anal. Appl.
  \textbf{18} (1967), no.~1, 75--84.

\bibitem{Gourdon2010}
X.~Gourdon, \emph{{The $10^{13}$ First Zeros of the Riemann Zeta Function, and
  Zeros Computation at Very Large Height}},
  \textbf{http://numbers.computation.free.fr/Constants/Miscellaneous/ zetazeros1e13-1e24.pdf}.

\bibitem{Helfgott2012}
H.A. Helfgott, \emph{{Minor arcs for Goldbach's problem}}, arXiv preprint
  arXiv:1205.5252 (2012).

\bibitem{Helfgott2013}
\bysame, \emph{{Major arcs for Goldbach's problem}}, arXiv preprint
  arXiv:1305.2897 (2013).

\bibitem{IEEE754}
IEEE, \emph{{IEEE Standard for Binary Floating-Point Arithmetic, IEEE Std
  754-1985.}}, 1985.

\bibitem{Lambov2008}
B.~Lambov, \emph{{Reliable Implementation of Real Number Algorithms: Theory and
  Practice}}, Lecture Notes in Computer Science, ch.~Interval Arithmetic Using
  SSE-2, Springer, 2008.

\bibitem{Moore1965}
R.E. Moore, \emph{Error in digital computation}, vol.~I, pp.~61--130, Wiley,
  1965.

\bibitem{Moore1966}
\bysame, \emph{{Interval analysis}}, vol.~60, Prentice-Hall Englewood Cliffs,
  New Jersey, 1966.

\bibitem{Muller2010}
J.M. Muller, \emph{{Correctly Rounded Mathematical Library}},
  \textbf{http://lipforge.ens-lyon.fr/www/ crlibm/}.

\bibitem{Odlyzko1988}
A.M. Odlyzko and A.~Sch{\"o}nhage, \emph{{Fast algorithms for multiple
  evaluations of the Riemann zeta function}}, Trans. Amer. Math. Soc.
  \textbf{309} (1988), no.~2, 797--809.

\bibitem{Rademacher1959}
H.~Rademacher, \emph{{On the
  Phragm$\acute{\textrm{e}}$n-Lindel$\ddot{\textrm{o}}$f theorem and some
  applications}}, Math. Z. \textbf{72} (1959), no.~1, 192--204.

\bibitem{Revol2002}
N.~Revol and F.~Rouillier, \emph{{A library for arbitrary precision interval
  arithmetic}}, 10th GAMM - IMACS International Symposium on Scientific
  Computing, Computer Arithmetic, and Validated Numerics, 2002.

\bibitem{Rumely1993}
R.~Rumely, \emph{{Numerical Computations Concerning the ERH}}, Math. Comp.
  \textbf{61} (1993), no.~203, 415--440.

\bibitem{Trudgian2011}
T.~Trudgian, \emph{{Improvements to Turing's method}}, Math. Comp \textbf{80}
  (2011), no.~276, 2259--2279.

\bibitem{Turing1953}
Alan~M. Turing, \emph{{Some calculations of the Riemann zeta-function.}}, Proc.
  Lond. Math. Soc. \textbf{3} (1953), no.~3, 99--117.

\bibitem{Walker1991}
J.S. Walker, \emph{{Fast Fourier Transforms}}, CRC press Boca Raton, 1991.

\end{thebibliography}
\end{document}